\documentclass[11pt,a4paper]{article}
\usepackage[cp866]{inputenc}
\usepackage[english]{babel}
\usepackage{amssymb,amsmath,amsfonts}
\usepackage{graphicx}

\date{}
\begin{document}

\title{On $\mathbb{Z}_4$-linear Reed-Muller like codes}
\author{Faina I. Solov'eva\thanks{This
research was partially supported by  the project N512 of the
Russian education program  "Development of scientific level of
Universities". }\\
Sobolev Institute of Mathematics\\
pr. ac. Koptyuga 4, Novosibirsk 630090\\
Russia\\
e-mail: sol@math.nsc.ru} \maketitle

\newtheorem {theorem}{Theorem}
\newtheorem {lemma}{Lemma}
\newtheorem {corol}{Corollary}
\newtheorem {pro}{Proposition}
\newcommand {\proofr }{{\par\medskip\noindent \bf Proof. }}
\newcommand {\proofend }{~~\raisebox{-1ex}{$\Box $}}
\newcommand {\point}{\hspace{-1.75mm}{\bf.\ }}
\newcommand {\pointt}{\hspace{-6.75mm}{\bf.\ }}

\begin{abstract}
For each $r$, $0\leq r \leq m$, it is presented the class of
quaternary linear codes $\mathcal{LRM}(r,m)$
 whose images under the Gray map
are binary codes with parameters of Reed-Muller $RM(r,m)$ code of order $r$.
\end{abstract}

\section{Preliminaries}

First consider some necessary definitions and notions. Consider
the ring $\mathbb{Z}_4$ of integers modulo $4$.  The set
$\mathbb{Z}_4^n$ is a module with addition operation  over the
ring $\mathbb{Z}_4$.
 The {\it Lee weight}, $w_L(\cdot)$,
of a quaternary vector is the sum of weights of its coordinate
positions:
 $$w_L(0)=0, w_L(1)=w_L(3)=1, w_L(2)=2.$$ The
{\it Lee distance}, $d_L(\cdot\,,\cdot)$, between any quaternary
vectors $x, y\in \mathbb{Z}_4^n$ is defined as $d_L(x,
y)=w_L(x-y)$. The set $\mathbb{Z}_4^n$ is a metric space with
respect to the  Lee metric.

A  {\it quaternary code} of length $n$ is a subset of the metric
space  $\mathbb{Z}_4^n$.  A quaternary code of length $n$ is {\it
linear} if it is a subgroup of the additive group of the
 ring $\mathbb{Z}_4^n$. We use further
capital letters for binary codes and  calligraphic for quaternary.

Let us remind the standard maps $\alpha$, $\beta$ and $\gamma$
from $\mathbb{Z}_4$ to $\mathbb{Z}_2$:
$$
\begin{array}[b]{c|ccc}
\mathbb{Z}_4 & \alpha & \beta & \gamma\\[1pt] \hline ~0~\rule{0pt}{12pt} & ~0~ & ~0~ & 0\\ 1
& 1 & 0 & 1\\ 2 & 0 & 1 & 1\\ 3 & 1 & 1 & 0
\end{array}.
$$
These maps can be extended in the usual way to maps from
$\mathbb{Z}_4^n$ to $\mathbb{Z}_2^n.$  The {\it Gray map}
$\phi\colon \mathbb{Z}_4^n \to \mathbb{Z}_2^{2n}$ is defined by
\[
\phi(x)=(\beta(x),\gamma(x)),\quad \mbox{for any}\ x\in
\mathbb{Z}_4^n.
\]
It is well known  that $\phi$ is an isometry of the metric spaces
 $\mathbb{Z}_4^n$ and $\mathbb{Z}_2^{2n}$.

  Two quaternary codes,
$\mathcal{C}$ and $\mathcal{D}$, of length $n$ are {\it
equivalent} if there exist a vector $x\in\mathbb{Z}_4^n$, a
permutation $\pi$ of $n$ coordinate positions and an inversion
$\sigma$ on $n$ coordinates such that $C=\pi(\sigma (D))+x$.  A
binary code is called {\it $\mathbb{Z}_4$-linear} with respect to
the Gray map if there exists an equivalent code $C$ such that its
preimage $\phi^{-1}(C)$ is linear.

 There are several known classes of
nonlinear binary codes with good properties which can be
represented as linear quaternary codes. Among them there are such
prominent codes as Preparata, Kerdock, Delsarte-Goethals,
Goethals-Delsarte,  some perfect codes, some Hadamard codes, see
the list of references
\cite{N1989,KN1992,HKCSS1994,NK1996,K2000,K2001,W1997}.

In \cite{K2000,K2001} the classifications of $\mathbb{Z}_4$-linear
perfect and $\mathbb{Z}_4$-linear  Hadamard codes are presented.
It is established that for any $n=2^k, n\geq 16$ the number of
nonequivalent $\mathbb{Z}_4$-linear perfect (Hadamard) codes is
$\lfloor(k+1)/2\rfloor$. All quaternary linear codes whose images
under the Gray map are
 perfect codes
can be described using Mollard construction, see  \cite{K2000}.

The representation of $\mathbb{Z}_4$-linear Preparata codes is
done in \cite{T2005}. Using switching approach it is established
that the set of all quaternary linear Preparata codes of length
$n=2^m$, $m$ odd, $m\ge 3$, is nothing more than the set of codes
of the form $\mathcal{H_{\lambda,\psi}}+\mathcal{M}$ with
$$
\mathcal{H_{\lambda, \psi}}=\{y+T_{\lambda}(y)+S_{\psi}(y)\mid
y\in H^n\},\quad \mathcal{M}=2H^n,
$$
where $T_{\lambda}(\cdot)$ and $S_{\psi}(\cdot)$ are vector fields
of a special form defined over the binary extended linear Hamming
code $H^n$ of length $n$. An upper bound on the number of
nonequivalent quaternary linear Preparata codes of length $n$ is
obtained, namely, $2^{n\log_2 n}$.

There are  several papers devoted to quaternary Reed-Muller
$RM(r,m)$ codes of order $r$, $0\leq r \leq m$. In
\cite{HKCSS1994} it is established that binary Reed-Muller
$RM(r,m)$ codes of order $r$, $r \in \{0,1,2,m-1, m\}$ are
$\mathbb{Z}_4$-linear and conjectured that all other Reed-Muller
codes are not $\mathbb{Z}_4$-linear. The conjecture is proved in
\cite{HLK1998}.

In \cite{HKCSS1994}  the class of quaternary codes
$\mathcal{QRM}(r,m)$ for each  $r$, $0\leq r \leq m$ is
introduced. The image of the code $\mathcal{QRM}(r,m)$ under the
map $\alpha$ (see the definition of the map above) is linear
Reed-Muller $RM(r,m)$ code for all $r$, $0\leq r \leq m$. The
class of the codes includes the quaternary linear Kerdock codes
and its dual the quaternary linear Preparata code from
\cite{HKCSS1994}. The generalization of the result is given in
\cite{BFP2005}. The class of the codes obtained in \cite{BFP2005}
includes all the quaternary linear Kerdock codes and the
quaternary linear Preparata codes. Thus,  another representation
of the quaternary linear Preparata codes is given. The images of
all these codes under the  map $\alpha$  are also linear
Reed-Muller codes.

In \cite{PR1997}  the additive Reed-Muller code
$\mathcal{ARM}(r,m)$ of order $r$, $0\leq r \leq m$  is defined.
The code is an additive subgroup of $\mathbb{Z}_2^{k_1}\times
\mathbb{Z}_4^{k_2},$ $k_1=2^{m-1},\, k_2=2^{m-2}, m\geq 2$ and it
is announced   that its binary image is linear for $r\in
\{0,1,m-1,m\}$ and nonlinear for $r=m-2, m > 3$.

\section{Reed-Muller-like codes}

In this section for every integer $r$, $r \in \{0,1,\ldots ,m\}$,
we construct the class of quaternary linear codes of length
$2^{m-1}$, code distance $d=2^{m-r}$ and size $2^k$, where
$$k = 1+ \left( \begin{array}{c} m \\
1 \end{array} \right) + \cdots + \left( \begin{array}{c} m \\
r \end{array} \right).$$ The image of any such code under the Gray
map is a binary  (not necessary linear) code with parameters of
Reed-Muller $RM(r,m)$ code of order $r$.

Let $v=(v_1,\ldots,v_m)$ range over $\mathbb{Z}_2^m$. The binary
{\it Reed-Muller code $RM(r,m)$ of order  $r$} is generated by all
binary vectors of length $2^m$ corresponding to the Boolean
functions $f(v)$ equaled to monomials of degree not more than $r$.
The code has the following parameters:

\noindent $\bullet$ length of the code is $n=2^m$;

\noindent $\bullet$ the size of the code is $2^k$, where   $k = 1+ \left( \begin{array}{c} m \\
1 \end{array} \right) + \cdots + \left( \begin{array}{c} m \\
r \end{array} \right);$

\noindent $\bullet$ the code distance $d=2^{m-r}.$

It is known, see \cite{MWSl},  that the binary  Reed-Muller code
$RM(r,m)$ of order  $r$ can be described by Plotkin (doubling)
construction:
$$RM(r,m) = \{(x,x+y) \mid x \in RM(r,m-1), y \in
 RM(r-1,m-1)\}.$$
%There are some switching and concatenation constructions for
% binary codes not necessary linear  with parameters of
%Reed-Muller $RM(r,m)$ codes, see, for example \cite{P1974,S}.

Binary not necessary linear code whose  parameters coincide with
parameters of the binary linear  Reed-Muller code $RM(r,m)$ of
order  $r$ we will call {\it Reed-Muller-like code} of order  $r$.
We are going to prove that among them there are
$\mathbb{Z}_4$-linear codes. Such codes have some regular
properties. For example, all such binary codes are transitive.
Preimages of them under the Gray map are quaternary codes with
parameters:

\noindent $\bullet$ length of the code is $n=2^{m-1}$;

\noindent $\bullet$ the size of the code is $2^k$, where $k = 1+ \left( \begin{array}{c} m \\
1 \end{array} \right) + \ldots +\left( \begin{array}{c} m \\
r \end{array} \right);$

\noindent $\bullet$ the code distance $d=2^{m-r}.$

We will denote any such quaternary linear code by
$\mathcal{LRM}(r,m)$, its binary image --  by $LRM(r,m)$.

We construct a sequence of quaternary linear $\mathcal{LRM}(r,m)$
codes that includes the following classes of quaternary linear
codes: the quaternary repetition code $\mathcal{LRM}(0,m)$, some
quaternary linear Hadamard codes from \cite{K2001}, some
quaternary linear extended perfect codes from \cite{K2000,K2001},
full-even weight code  $\mathcal{LRM}(m-1,m)$ and
$\mathcal{LRM}(m,m)=\mathbb{Z}_4^{2^m}$. We construct the
$\mathcal{LRM}(r,m)$ codes  by induction on $m$, where $m=\log n,
m>1$. For $m=1$ there exist the following trivial quaternary
linear codes:

a) the quaternary linear code
$\mathcal{LRM}(0,1)=\mathcal{RM}(0,1)=\{(0),(2)\},$ its binary
image is
$LRM(0,1)=\phi(\mathcal{LRM}(0,1))=RM(0,1)=\{(0,0),(1,1)\};$

b) the quaternary linear code
$\mathcal{LRM}(1,1)=\mathcal{RM}(1,1)=\{(0),(1),(2),(3)\},$ its
binary image is
$LRM(1,1)=\phi(\mathcal{LRM}(1,1))=RM(1,1)=\{(0,0),(0,1),(1,1),((1,0)\}.$

Let us consider also more interesting case $m=2$. Here we have the
following quaternary linear codes and their binary images under the Gray
map:

 a) the quaternary linear code
$\mathcal{LRM}(0,2)=\mathcal{RM}(0,2)=\{(0,0),(2,2)\}$ and its
binary image
$LRM(0,2)=\phi(\mathcal{LRM}(0,2))=RM(0,2)=\{(0,0,0,0),(1,1,1,1)\};$

 b) the quaternary linear code
$$\mathcal{LRM}(1,2)=\mathcal{RM}(0,2)=\{(0,0),(1,1),(2,2),(3,3),(1,3),(3,1),(0,2),(2,0)\}$$
and its binary image full-even weight  code of length 4:
$LRM(1,2)=\phi(\mathcal{LRM}(1,2))=RM(1,2);$

 c) the quaternary linear code
$\mathcal{LRM}(2,2)=\mathbb{Z}_4^2$ and its binary image
$$LRM(2,2)=\phi(\mathcal{LRM}(2,2))=RM(2,2).$$

Let $\mathcal{LRM}(r,m-1)$ and  $\mathcal{LRM}(r-1,m-1)$ be any
two quaternary linear codes with parameters
$$(n=2^{m-2}, \, 2^k, \, d=2^{m-r-1})\,\,\, \mbox{and} \,\,\,
(n=2^{m-2}, \, 2^{k^{\prime}}, \, d=2^{m-r}),$$ where $$k=\sum_{i=0}^{r} \left( \begin{array}{c} m-1 \\
i \end{array} \right), k^{\prime}=\sum_{i=0}^{r-1} \left( \begin{array}{c} m-1 \\
i \end{array} \right).$$ It is not difficult to show that Plotkin
construction applied to these codes gives us a quaternary linear
$\mathcal{LRM}(r,m)$ code of order $r$, i.e.
$$\mathcal{LRM}(r,m)= \{(x,x+y) \mid x \in \mathcal{LRM}(r,m-1), y \in
 \mathcal{LRM}(r-1,m-1)\}.$$
 So we get the following result.

\begin{theorem}\point For any $r$, $0\leq r \leq m$, $m\geq 1$,  the
set $\mathcal{LRM}(r,m)$ is a quaternary linear code with
parameters
\begin{equation}(n=2^{m-1}, \, 2^k, \, d=2^{m-r}),\, \mbox{where} \,\, k=\sum_{i=0}^{r} \left( \begin{array}{c} m \\
i \end{array} \right), \end{equation}
 whose image under the Gray map is a binary
code with parameters of Reed-Muller $RM(r,m)$ code of order $r$.
\end{theorem}

Taking into account that there exist quaternary linear Hadamard
codes, all quaternary linear extended perfect codes from
\cite{K2000,K2001} whose binary images under the Gray map  are
nonlinear codes and all $RM(r,m)$ codes of order $r$, $r \in
\{3,\ldots, m-2\}$ are not $\mathbb{Z}_4$-linear
 we get a sequence of quaternary linear codes $\mathcal{LRM}(r,m)$ of
order $r$ such that all their binary images are nonlinear codes
with parameters of Reed-Muller codes $RM(r,m)$  for any order $r\in
\{3,\ldots, m-2\}$.

Let us show that the codes from Theorem 1 do not equivalent to
codes from \cite{BFP2005} for $r\in \{3,\ldots, m-2\}$. According
to \cite{BFP2005} a quaternary linear code $\mathcal{QRM}(r,m-1)$
of length $2^{m-1}$ has size $2^{2k}$, where $k=\sum_{i=0}^{r} \left( \begin{array}{c} m - 1\\
i \end{array} \right)$. From this fact and  (1) we conclude that any
quaternary linear codes $\mathcal{LRM}(r,m)$ and
$\mathcal{QRM}(r,m)$ of the same length have different sizes, so
they do not equivalent to each other. All binary images of the
codes $\mathcal{LRM}(r,m)$ under the Gray map are nonlinear  codes
with parameters of Reed-Muller codes $RM(r,m)$ of order $r\in
\{3,\ldots, m-2\}$, but all binary images of codes
$\mathcal{QRM}(r,m)$ are linear $RM(r,m)$ codes of the same order.
Therefore their binary images having the same parameters are also
nonequivalent to each other.

\end{document}